\documentclass{article}
\title{On length and product of harmonic forms in K\"ahler geometry
\footnote{
\newline
This research was supported by the VW-Research Group "Special geometries in Mathematical Physics " and EDGE, Research Training Network HPRN-CT-2000-00101, supported by the European Human Potential Programme}}
\author{Paul-Andi Nagy }
\date{\today}
\usepackage{amsfonts,amssymb}
\usepackage[backref]{hyperref}

\newtheorem{teo}{Theorem}[section]
\newtheorem{lema}{Lemma}[section]
\newtheorem{pro}{Proposition}[section]
\newtheorem{defi}{Definition}[section]
\newtheorem{rema}{Remark}[section]
\newtheorem{coro}{Corollary}[section]
\newtheorem{nr}{}[section]
\newtheorem{ex}{Example}[section]
\begin{document}
\maketitle
\abstract{ \normalsize
Motivated by understanding the limiting case of a certain systolic inequality we study compact Riemannian manifolds having all harmonic $1$-forms of constant
length. We give complete characterizations as far as K\"ahler and hyperbolic geometries are concerned. In the second part of the paper, we give algebraic and topological obstructions to the existence
of a geometrically $2$-formal K\"ahler metric, at the level of the second cohomology group. A strong interaction with almost K\"ahler geometry is to be noted. In complex dimension $3$, we list all the possible values of the second Betti number of a geometrically $2$-formal K\"ahler metric.
\large
\tableofcontents
\section{Introduction}
Let $(M^n,g)$ be a compact oriented Riemannian manifold. We denote by $\Lambda^{p}(M), 0 \le  p \le n$ the space of smooth, real valued, $p$-forms of
$M$. The standard deRham complex
$$ \ldots \rightarrow \Lambda^{p}(M) \stackrel{d}{\rightarrow}  \Lambda^{p+1}(M) \rightarrow \ldots $$
where $d$ stands for exterior derivative is then used to introduce the deRham cohomology groups, to be denoted by $H_{DR}^{p}(M)$. The topological
information contained in these cohomology groups may be understood geometrically, using Hodge theory, by means of the isomorphisms
\begin{nr} \hfill
$ H_{DR}^{p}(M) \equiv {\cal{H}}^p(M,g), \ 0 \le p \le n. \hfill $
\end{nr}
Here the space of {\it{harmonic}} $p$-forms of $(M^n,g)$ is defined by
$$ {\cal{H}}^p(M,g)=\{ \alpha \in \Lambda^p(M) : \Delta \alpha=0 \}. $$
The Laplacian on forms is given by $\Delta=dd^{\star}+d^{\star}d$ where $d^{\star}$ is the formal adjoint of $d$ with respect to the
given metric and orientation of $M$.\par
In this paper we investigate various notions of "constancy" related to harmonic forms. The first one is introduced by the following
\begin{defi}
Let $(M^n,g)$ be compact and oriented. It is said to satisfy the hypothesis $(CL_p)$ for some $1 \le p \le n-1$ iff
every harmonic $p$-form has pointwisely constant norm.
\end{defi}
Manifolds satisfying hypothesis $(CL_1)$ appear to be naturally related to a generalized systolic inequality. More precisely,
for a compact, orientable Riemannian manifold $(N^n,g)$ with
non-vanishing first Betti number one defines the stable $1$-systole $stsys_1(g)$ in terms of the stable norm (see \cite{Ban1,Ban2} for details). Let
$sys_{n-1}(g)$ be the infimum of the $(n-1)$-volumes of all nonseparating hypersurfaces in $N$. Then
the following systolic inequality, previously established in \cite{Hebda} in the case when the first Betti
number equals $1$ holds (see \cite{Ban1}) :
\begin{nr} \hfill \label{sys}
$ stsys_1(g)\cdot sys_{n-1}(g) \le \gamma^{\prime}_{b_1(N)}\cdot vol(g).\hfill $
\end{nr}
Here $\gamma^{\prime}_{b_1(N)}$ is the Berg\'e-Martinet constant for whose definition we send again the reader to \cite{Ban2}. The important point for
us is that it was shown in \cite{Ban2} that if equality in (1.2) occurs then $(N^n,g)$ satisfies the hypothesis $(CL_1)$. Note that the converse is false, as flat tori always satisfy
$(CL_1)$ but saturate (\ref{sys}) iff they are dual-critical.\par
Riemannian manifolds $(N^n,g)$ saturating (\ref{sys}) have strong geometric properties. It was proved in \cite{Ban2}, Thm. 1.2, that in this case $(N^n,g)$ is the total space of a Riemannian submersion with minimal fibers to a flat torus, whose projection is actually the Albanese map.
Therefore, in the special case when $b_1(N)=n-1$ it follows that the fibers of the Albanese map must be totally geodesic. Using Chern-Weil theory and an argument that reproduces in part that in section 6 of \cite{Ban2}, we showed in \cite{NV} that the only possible topologies of manifolds $N^n$ which admit a metric satisfying $(CL_1)$ and have $b_1(N)=n-1$ are those of $2$-step nilmanifolds with $1$-dimensional kernel. Equivalently, the above class of manifolds is parametrized by couples
$(T, \omega)$ where $T$ is a flat $(n-1)$-torus and $\omega$ is a non zero, integral cohomology class on $T$.
\par
For a compact oriented Riemannian manifold $(M,g)$ we now set
$$ H^{\star}_{DR}(M)=\bigoplus \limits_{p \ge 0}H_{DR}^p(M) \ \mbox{and} \ {\cal{H}}^{\star}(M,g)= \bigoplus \limits_{p \ge 0} {\cal{H}}^{p}(M,g)$$
Whilst $H^{\star}(M)$ is a graded algebra, in general ${\cal{H}}^{\star}(M,g)$ is not an algebra with respect to the
wedge product operation for there is no reason the isomorphism (1.1) descends to the level of harmonic forms. Our next definition
is related to this fact.
\begin{defi}
Let $(M^n,g)$ be compact and oriented. \\
(i) The metric $g$ is $p$-formal for some $1 \le p \le n-1$ iff the product of any
harmonic $p$-forms remains harmonic. \\
(ii) The metric $g$ is formal iff the product of any two harmonic  forms remains harmonic.
\end{defi}
Following \cite{Kot} we also set
\begin{defi} Let $M^n$ be compact and oriented. $M^n$ is geometrically formal iff it admits a formal Riemannian metric.
\end{defi}
From a topological viewpoint, geometric formality implies that the rational
homotopy type of the manifold is a formal consequence of the cohomology ring \cite{Sul}. Basic examples
are compact Riemannian symmetric spaces. In fact, in the recent \cite{Kot}, it was proved that in dimension
$3$ and $4$ every geometrically formal manifold has the real cohomology algebra of a compact Riemannian symmetric
space. In higher dimensions, there are very few
general facts known about geometrically formal manifolds; for instance
formal metrics satisfy hypothesis $(CL_p)$ for all $1 \le p \le n-1$ \cite{Kot}. Note that, by contrast, the class of
(non necessarily invariant) metrics on nilmanifolds studied in \cite{NV} satisfy hypothesis $(CL_p)$ whenever $1 \le p \le n-1$ but none of the $p$-formality hypothesis.
Moreover, it is known that certain classes of homogeneous spaces fail
to be geometrically formal for cohomological reasons \cite{Kot1}. \par
In this note we place ourselves in the context of K\"ahler manifolds and we investigate geometric consequences of the constant length hypothesis
and of (geometric) formality for low degree harmonic forms. Our paper is organized as follows. In section 2 we prove the following.
\begin{teo}
Let $(M^{2n},g,J)$ be a compact K\"ahler manifold. Then every harmonic $1$-form of pointwisely constant length is parallel with respect to the Levi-Civita connection 
of $g$. In particular, if $g$ satisfies the hypothesis $(CL_1)$ then $(M^{2n},g,J)$ is locally the 
Riemannian (and biholomorphic) product of a compact, simply connected K\"ahler manifold and of a flat torus.
\end{teo}
Note that the result of theorem 1.1 is no longer available if instead of having a K\"ahler structure we require only the presence of an almost K\"ahler
one (see section 2 for an example). It also follows that a compact K\"ahler manifold which is locally irreducible and not flat never
saturates the systolic inequality (1.2). Moreover in section 2 we remark that the length of a harmonic $1$-form on a compact hyperbolic manifold 
cannot be constant; this is actually a consequence of a result in \cite{Kim} and holds in fact for compact locally symmetric spaces of negative curvature 
\cite{Wal}. We propose a different, very simple proof. \par
\par
The rest of the paper is concerned with the study of obstructions to geometric $2$-formality. Note however that
every compact K\"ahler is {\it{topologically}} formal by results in \cite{Deligne}. In section 3, we show that harmonic $2$-forms of a
$2$- formal K\"ahler manifold have a global spectral decomposition and constant eigenvalues. This is enforcing the
opinion, already presented in \cite{Kot} that geometric formality is weakening the notion of Riemannian holonomy reduction. Based upon this we are able to prove the following.
\begin{teo}
Let $(M^{2n},g,J)$ be a compact K\"ahler manifold and assume that the metric $g$ is $2$-formal. Then : \\
(i) the space ${\cal{H}}^{1,1}$ of $J$-invariant harmonic forms is spanned by almost K\"ahler forms
compatible with the metric $g$. \\
(ii) the space ${\cal{H}}^2_{-}$ of $J$-anti-invariant harmonic two-forms consists only in parallel forms.
\end{teo}
By contrast, recall that simply connected, irreducible compact Hermitian symmetric space have second Betti number
equal to $1$.
As an application of theorems 1.1 and 1.2 we show in section 4 that
geometrically formal K\"ahler manifold having a maximal Betti number are flat. 
Further consequences of theorem 1.2 are investigated under various curvature
assumptions in section 4. For example, we prove that locally irreducible
$8$-dimensional hyperk\"ahler manifolds cannot be geometrically formal. In
section 5 we study the case  of geometrically formal K\"ahler
manifolds of complex dimension $3$. We are able to give the possible values of
the second Betti number in this situation together with more precisions
concerning the algebraic structure of the second cohomology group. We prove :
\begin{teo} Let $(M^6,g,J)$ be a geometrically formal K\"ahler manifold. If the metric $g$ is 
locally irreducible then $b_1(M)=b_2^{-}(M)=0$. Moreover one has $b_2(M) \le 3$ and ${\cal{H}}^{1,1}$ is
spanned by mutually commuting almost K\"ahler structures.   
\end{teo} 
As a final remark, we mention that finding further, first order obstructions
to geometric formality in the K\"ahler case relies on understanding
the algebraic structure of the space of harmonic $p$-forms, $p\ge 3$.
\section{The length of harmonic $1$-forms}
This section will be devoted to the investigation of geometric issues of the existence of a harmonic
$1$-form of constant length on a compact K\"ahler manifold. In particular our discussion 
will lead to the proof of Theorem 1.1. Before proceeding we need to recall some basic 
material related to a particular class of foliations. \par
Let $(M,g)$ be a Riemannian manifold equipped with a smooth foliation ${\cal{F}}$ and let us denote 
by ${\cal{V}}$ the integrable distribution on $M$ induced by ${\cal{V}}$. We consider the splitting 
\begin{nr} \hfill
$ TM={\cal{V}} \oplus H \hfill $
\end{nr}
where $H$ is the orthogonal complement of ${\cal{V}}$. From now on we will denote 
by $V,W$ vector fields in ${\cal{V}}$ and by $X,Y,Z$ etc. vector fields in $H$. Let $\nabla$ be 
the Levi-Civita connection of the metric $g$. Recall that $H$ is {\it{totally geodesic}} iff 
$\nabla_XY$ belongs to $H$. 
Foliations ${\cal{F}}$ satisfying this condition -to be assumed, unless otherwise stated, in the rest of our preliminaries-
shall be termed {\it{transversally totally geodesic}}.
Then we note that ${\cal{F}}$ is a particular kind of {\it{Riemannian foliation}}, meaning that 
$$ ({\cal{L}}_Vg)(X,Y)=0.$$
We will present below some basic notions related to this class of foliations, following closely \cite{Besse,Tond}. To 
begin with, let $\overline{\nabla}$ be the 
orthogonal projection of $\nabla$ onto the splitting (2.1). Then it is easy to verify that 
$\overline{\nabla}$ defines a metric connection (with torsion) preserving the distributions ${\cal{V}}$ and 
$H$. \par
An important object in our study will be the O'Neill tensor $T$ defined by (see \cite{Tond}, p.49)
$$ T_EF=(\nabla_{E_{\cal{V}}}{F_{\cal{V}}})_H+(\nabla_{E_{\cal{V}}}F_H)_{{\cal{V}}}
$$
whenever $E,F$ belong to $TM$; here the subscript denotes orthogonal projection on the subspace. It follows 
that $T$ vanishes on $H \times H$ and $H \times {\cal{V}}$, it is symmetric on ${\cal{V}} \times {\cal{V}}$ (since 
${\cal{V}}$ is integrable) and furthermore we have $<T_VX,W>=-<X,T_VW>$. \par
Based on these definitions it is easy to check that the connections $\nabla$ and $\overline{\nabla}$ are related 
to the tensor $T$ by : 
$$\begin{array}{ll}
\nabla_XY=\overline{\nabla}_XY & \nabla_XV=\overline{\nabla}_XV \vspace{2mm} \\
\nabla_VW=\overline{\nabla}_VW+T_VW & \nabla_VX=\overline{\nabla}_VX+T_VX.
\end{array} $$
The last notion needed for our purposes is related to the curvature $\overline{R}$ of the connection 
$\overline{\nabla}$ defined by $\overline{R}(E,F)=\overline{\nabla}_{[E,F]}-\overline{\nabla}_E \overline{\nabla}_F+
\overline{\nabla}_F \overline{\nabla}_E$ for all vector fields $E$ and $F$ of $M$. Then the {\it{transversal Ricci tensor}} 
$Ric^H :H \to H$ is given by 
$$ <Ric^HX,Y>=\sum \limits_{e_i \in H}^{} \overline{R}(X,e_i,Y,e_i)$$
for an arbitrary local orthonormal basis $\{ e_i \}$ of $H$. \par
$\\$
We consider now a compact K\"ahler manifold $(M^{2n},g,J)$ admitting a harmonic $1$-form $\alpha$ of pointwisely constant 
length. Let $\zeta$ be the vector field dual to $\alpha$ and consider the distribution $H$ spanned by $\zeta$ and $J\zeta$. Moreover, 
let ${\cal{V}}$ be the orthogonal complement of $H$ in $TM$. Our starting point toward the proof of Theorem 1.1 is the following 
\begin{lema} The distribution ${\cal{V}}$ is integrable and the distribution $H$ is totally geodesic. Moreover we have $Ric^H=0$.
\end{lema} 
{\bf{Proof}} : \\ 
If $\beta$ is a $1$-form on $M$ we let $J$ act on $\beta$ by $(J\beta)X=\beta(JX)$ for all $X$ in $TM$. Since $M$ is compact 
we know that $J\alpha$ must be closed. Together with the closedeness of $\alpha$ this leads to the integrability of ${\cal{V}}$. 
By construction, we have that the splitting $TM={\cal{V}} \oplus H$ is $J$-invariant and moreover on a compact K\"ahler 
manifold any harmonic $1$-form $\gamma$ is holomorphic, that is 
\begin{nr} \hfill 
$ \nabla_{JX}\gamma=J\nabla_X \gamma\hfill $
\end{nr}
whenever $X$ belongs to $TM$, where $\nabla$ is the Levi-Civita of the metric $g$. Therefore $H$ is a holomorphic distribution and we use results in
\cite{Wat}, to conclude that $H$ is totally geodesic. Furthermore $H$ is actually transversally flat; this follows easily from (2.2) and the fact that 
$\alpha$ has constant length and leads to the last assertion of our Lemma.
$\blacksquare$ \par

But the geometry of the foliations satisfying the conditions in Lemma 2.1 can be completely 
ruled out. In fact we shall prove the slightly more general 
\begin{pro} Let $(M^{2n},g,J)$ be a compact K\"ahler manifold, supporting a
foliation with complex leaves which is transversally totally geodesic  and
with non-negative transverse Ricci curvature. Then $M$ is locally a Riemannian
(and K\"ahler ) product. 
\end{pro} 
{\bf{Proof}} : \\ 
Let ${\cal{V}}$ be the distribution tangent to the leaves of the foliation and $H$ its orthogonal 
complement. The splitting $TM={\cal{V}} \oplus H$ is then orthogonal and $J$-invariant. Let $\nabla$ be the 
Levi-Civita connection of the metric $g$. Because $(g,J)$ 
is K\"ahler the connection $\overline{\nabla}$ is Hermitian, i.e $\overline{\nabla}J=0$ and also the O'Neill 
tensor $T$ of the foliation satisfies $[T_V,J]=0$. Since $H$ is a totally
geodesic distribution we have \cite{Besse} : 
$$<(\overline{\nabla}_XT)(V,W),Y>=<(\overline{\nabla}_YT)(V,W),X>.$$
It follows that
$(\overline{\nabla}_{JX} T)(JV,W)=(\overline{\nabla}_XT)(V,W)$.
Derivating in the direction of $e_i$ (here $\{ e_i \}$ is an arbitrary local
orthonormal basis in $H$) we get : 
$$ (\nabla_H^{\star} \nabla_H)T=\frac{1}{2}
J \sum \limits_{e_i \in H}^{}\overline{R}(e_i, Je_i).T$$ where $\nabla_H$
denotes derivation with respect to $\overline{\nabla}$, in the direction of
$H$.  The tensor $\overline{R}$, the curvature tensor of the connection
$\overline{\nabla}$, acts on $T$ by
$$(\overline{R}(X,Y).T)(V,W)=\overline{R}(X,Y)[T_VW]-T_{\overline{R}(X,Y)V}W-T_V \overline{R}(X,Y)W.$$ To compute this last term we note that (see \cite{Tond})
$$\overline{R}(X,Y,V,W)=R(X,Y,V,W)=<T_WX,T_VY>-<T_VX,T_WY>$$ (see \cite{Tond}). After a short computation this yields to  
$$ \frac{1}{2} J \sum \limits_{e_i \in H}^{} (\overline{R}(e_i, Je_i).T)(V,W)=-(Ric^H(T_VW)+T_{SV}W +T_{V}SW) $$ where 
the symmetric endomorphism $S : {\cal{V}} \to {\cal{V}}$ is given by 
$<SV,W>= \sum \limits_{e_i \in H}^{} <T_Ve_i, T_We_i>$. Taking 
the scalar product with $T$ implies by means of the positivity of the
transversal Ricci curvature that $<(\nabla_H^{\star} \nabla_H)T,T> \le -2\vert
S \vert^2$. The vanishing of $T$ (and hence the $\nabla$-parallelism of the splitting $TM={\cal{V}} \oplus H$) follows  now simply by integration over $M$
followed by a positivity argument. $\blacksquare$ 
\begin{rema}
(i) Proposition 2.1 actually holds when relaxing the hypothesis on the foliation to transversal integrability. Since the proof 
is more involved and not directly related to our present investigations we chose not to present it here.\\
(ii) A result similar to Proposition 2.1 was proved in \cite{Apo1}, under the assumption that $Ric$ has constant positive 
eigenvalues on ${\cal{V}}$ and $H$, by making use of Sekigawa's integral formula.
\end{rema}
The following example shows that the splitting result in Theorem 1.1 is intimately related to the presence of a K\"ahler structure and cannot hold in presence
of an almost K\"ahler, non-K\"ahler, structure.
\begin{ex}
Let $(N,g)$ be a $3$-dimensional geometrically formal manifold. Many non-symmetric examples are known to exist \cite{Kot}, and in particular we must have $b_1(N)=1$. Let $\alpha$ be the harmonic $1$-form of $N$  whose
length equals $1$ and let $M=\mathbb{S} ^1 \times N$ be endowed with the product metric. It is a simple verification
that $\omega=dt \wedge \alpha +\star \alpha$ defines a compatible symplectic form giving $M$ the structure of a compact almost K\"ahler manifold. But $b_1(M)=2$ and in general $N$ is a locally irreducible Riemannian manifold.
\end{ex}
Higher dimensional examples, endowed with non-homogeneous Riemannian metrics can be obtained by taking 
the product with $\mathbb{S}^1$ of the class of nilmanifolds studied in \cite{NV}.
\par
Combining Proposition 2.1 with Lemma 2.1 we are lead directly to the proof of Theorem 1.1.
As a direct consequence we obtain : 
\begin{coro}
Let $(M^{2n},g,J)$ be a locally irreducible K\"ahler manifold with $b_1(N)>0$. Then inequality (1.2) is always strict.
\end{coro}
We equally have :
\begin{coro}
Let $(M^{2n},g,J)$ be a compact K\"ahler manifold. If $g$ is $1$-formal then any harmonic $1$-form is parallel for the 
Levi-Civita connection of $g$.
\end{coro}
{\bf{Proof}} : \\
Let $\alpha$ be a harmonic $1$-form, with dual vector field $\zeta$. Because $\alpha$ and $J\alpha$ are co-closed a simple computation 
shows that the vector field dual to $d^{\star}(\alpha \wedge J\alpha)$ equals $[\zeta,J\zeta]$. Since $g$ is $1$-formal it follows that $[\zeta, J\zeta]=0$ and 
since $\alpha$ is holomorphic (i.e. it satisfies (2.2)) we arrive at $\nabla_{\zeta}\zeta=0$. The closedeness of $\alpha$ implies that 
$\alpha$ is of constant length hence the proof is completed by applying Theorem 1.1.
$\blacksquare$ \par
$\\$
The Riemannian submersion technique used below can be also used to disqualify some other locally symmetric spaces from having harmonic $1$-forms
of constant length. Also the following proposition happens to provide the answer to an open question in \cite{Ban2} as well as it provides a serious
obstruction to the geometric formality of compact hyperbolic manifolds. Note that the result below follows in fact directly from the more general
result in \cite{Kim}, asserting the non-existence of Riemannian submersions from compact hyperbolic manifolds. We give the (different) proof 
mainly because of its simplicity.
\begin{pro}
Let $(M^{2n},g), n\ge 1$ be a compact manifold with constant negative sectional curvature. Then any harmonic $1$-form of constant length vanishes.
\end{pro}
{\bf{Proof}} : \\
Let us suppose that $(M,g)$ admits a non-vanishing harmonic $1$-form $\alpha$ of constant length. Let $\zeta$ be the vector field dual to $\alpha$ and
assume, for simplicity, that $\zeta$ is of unit length. We define now $H$ to be the $1$-dimensional distribution spanned by $\zeta$ and let ${\cal{V}}$ be its orthogonal
complement in $TM$. \par
Let $\nabla$ be the Levi-Civita connection of the metric $g$. The compactness of $M$  imply that $\alpha$ is closed and co-closed. These two
equations imply easily (see \cite{Ban2} for
a related discussion) that
$\nabla_{\zeta} \zeta=0$ ($H$ is totally geodesic) and furthermore that ${\cal{V}}$ is integrable and also minimal. For the minimality of ${\cal{V}}$ will be quite important
for us we note that it is equivalent with the fact that $\alpha$ is coclosed. In other words, the splitting $TM={\cal{V}} \oplus H$ defines a 
transversally totally geodesic foliation (hence Riemannian) with 
leaves on $M$ and we are going to use O'Neill's structure equations for such an object. \par
Since $H$ is $1$-dimensional the O'Neill tensor $T$ can be written as 
$$T_VW=<SV,W>\zeta$$ 
where $S : {\cal{V}} \to {\cal{V}}$ is a symmetric and traceless tensor (because of the integrability and minimality of ${\cal{V}}$). If $R$ denotes 
the curvature tensor of the Levi-Civita connection, we recall that the following equation holds (see \cite{Tond}):
\begin{nr} \hfill
$ R(\zeta, V, \zeta, W)=-<(\overline{\nabla}_{\zeta}T)(V,W), \zeta>+<T_V{\zeta}, T_W{\zeta}>\hfill $
\end{nr}
whenever $V,W$ belong to ${\cal{V}}$, where $<T_{V}{\zeta}, W>=-<\zeta, T_VW>$. Taking into account that, after a suitable renormalization, we can assume that
$$ -R(X,Y,Z,U)=<X,U><Y,Z>-<X,Z><Y,U>$$
for all $X,Y,Z,U$ in $TM$, equation (2.3) can further rewritten as
\begin{nr} \hfill
$ <V,W>=-<(\overline{\nabla}_{\zeta}S)V,W>+<SV,SW>\hfill $
\end{nr}
for all $V,W$ in ${\cal{V}}$. Starting from $Tr(S)=0$, an elementary manipulation of (2.4) yields by induction to
$Tr(S^{2k+1})=0$ and $Tr(S^{2k})=n-1$ for all natural $k$. But the last relation implies immediately that $S^2=1_H$ and because 
$S$ is traceless we find that $n-1$ is even, a contradiction.
$\blacksquare$
\begin{coro}
Let $(M^{2n},g)$ be compact with constant negative sectional curvature. The following hold  : \\
(i) the inequality (1.2) is a strict one.\\
(ii) if $g$ is a formal metric we must have $b_1(M)=0$.
\end{coro}
We finish this section by pointing out the important fact
that both results of Proposition 2.2 and Corollary 2.3 hold in the more general context of compact locally symmetric spaces of negative (sectional) curvature
in virtue of results in \cite{Wal}. 

\section{Algebraic obstructions} 
In this section we are going to examine some elementary algebraic
obstructions to the existence of a $2$-formal K\"ahler metric. We begin by a brief review of some facts of K\"ahler geometry,
of relevance for our purposes. \par
For any compact K\"ahler manifold
$(M^{2n},g,J)$ we can consider the decomposition 
$$\Lambda^2(M)=\Lambda^{1,1}(M) \oplus \Lambda^2_{-}(M)$$ 
where $\Lambda^{2}_{-}(M)=\{
\alpha : J\alpha=-\alpha \}$. Here $J$ acts on a two form $\alpha$ by
$(J\alpha)(X,Y)=\alpha(JX,JY)$ whenever $X,Y$ belong to $TM$. The non-standard
notation is motivated by the fact that we are working with real-valued
differential forms. \par 
We have a further decomposition 
$$\Lambda^{1,1}(M)=\Lambda^{1,1}_{0}(M) \oplus C^{\infty}(M).\omega$$ 
where $\omega=g(J \cdot,
\cdot)$ is the K\"ahler form of $(g,J)$ and $\Lambda^{1,1}_{0}(M)$ is the
sub-bundle of $\Lambda^{1,1}(M)$ consisting of primitive forms. We denote now
by ${\cal{H}}^p, p \ge 0$ the space of harmonic $p$-forms with respect to the
metric $g$. The previous decompositions have analogues at the level of
harmonic forms 
\begin{nr} \hfill ${\cal{H}}^2={\cal{H}}^{1,1} \oplus
{\cal{H}}^{2}_{-} \hfill $ \end{nr} 
and 
\begin{nr} \hfill
${\cal{H}}^{1,1}={\cal{H}}^{1,1}_{0} \oplus \mathbb{R}.\omega \hfill $
\end{nr} 
with the obvious notational conventions. Moreover, we will denote by
$h^{1,1}$ the dimension of ${\cal{H}}^{1,1}$ and by  $b_2^{-}$ that of
${\cal{H}}^{2}_{-}$. \par Let ${\cal{S}}$ the space of  symmetric and
$J$-invariant endomorphisms $S$ of $TM$ having the property that  $<SJ \cdot,
\cdot>$ belongs to ${\cal{H}}^{1,1}$. Clearly, ${\cal{S}}$ and
${\cal{H}}^{1,1}$ are isomorphic and its worthwhile to note that all elements
of ${\cal{S}}$ have constant trace, in virtue of (3.2). In the same  way we
define the space ${\cal{A}}$ as the space of skew-symmetric, $J$-anti-commuting
endomorphisms of $TM$ which are associated to an element of
${\cal{H}}^{2}_{-}$. Note the important fact that $J \cdot {\cal{A}} \subseteq
{\cal{A}}$. \par Another aspect of K\"ahler geometry, of particular
significance for us, is that the operator $L$ defined as exterior
multiplication with the K\"ahler form preserves the space of the harmonic
forms of the manifold. This is a consequence of the fact that $L$ commutes
with the Laplacian acting on forms (see \cite{Gold}). Since the Laplacian is a
self-adjoint operator it also follows that $L^{\star}$, the adjoint of $L$,
preserves the space of harmonic forms. \par 
We now give a first set of elementary algebraic obstructions to the presence of 
a $2$-formal K\"ahler metric. If $A$ and $B$ are endomorphisms
of some vector bundle over a manifold we will denote by $\{A,B \}=AB+BA$ their
anti-commutator. The whole discussion in this section will be based on the
lemma below. 
\begin{lema} Let $(M^{2n},g,J)$ be a K\"ahler manifold such that the metric $g$ is 
$2$-formal. The following hold  : \\ 
\begin{nr} \hfill  $\{{\cal{S}},
{\cal{S}}\} \subseteq {\cal{S}}. \hfill $ 
\end{nr} and  
\begin{nr} \hfill
$\{{\cal{A}}, {\cal{A}}\} \subseteq {\cal{S}}. \hfill $ 
\end{nr} 
We also have 
\begin{nr} \hfill
$\{{\cal{S}}, {\cal{A}}\} \subseteq {\cal{A}}. \hfill $
\end{nr}
\end{lema}
{\bf{Proof}} : \\
Let us prove (3.3). Consider $\alpha$ and $\beta$ in ${\cal{H}}^{1,1}$ with associated symmetrics $S_1$ and $S_2$. We fix
$\{ e_i, 1 \le i \le 2n \}$ be a local orthonormal basis and write : 
$$ L^{\star}(\alpha \wedge \beta)=\frac{1}{2} \sum \limits_{i=1}^{2n}Je_i \lrcorner (e_i \lrcorner (\alpha \wedge \beta)).$$
But 
$$ 
\begin{array}{lr}
Je_i \lrcorner (e_i \lrcorner (\alpha \wedge \beta))=Je_i \lrcorner ((e_i \lrcorner \alpha )\wedge \beta+
\alpha \wedge (e_i \lrcorner \beta))= \vspace{2mm} \\
(Je_i \lrcorner e_i \lrcorner \alpha ) \cdot \beta -(e_i \lrcorner \alpha ) \wedge
(Je_i \lrcorner \beta ) +(Je_i \lrcorner \alpha ) \wedge (e_i \lrcorner \beta )+\alpha \cdot (Je_i \lrcorner e_i \lrcorner 
\beta).
\end{array} $$ 
Assuming the basis to be Hermitian we get $L^{\star} (\alpha \wedge \beta)=L^{\star}\alpha \cdot \beta+
\alpha \cdot L^{\star} \beta-\gamma$ where 
$$ \gamma=\sum \limits_{i=1}^{2n} (e_i \lrcorner \alpha ) \wedge (Je_i \lrcorner \beta).$$
Now a short computation shows that $\gamma=<\{S_1, S_2 \}J \cdot, \cdot >$ 
and since $L^{\star}(\alpha \wedge \beta) $ belongs to ${\cal{H}}^{1,1}$ whilst $L^{\star}\alpha, L^{\star}\beta$ are constants
we get that $\gamma$ is equally in ${\cal{H}}^{1,1}$ and the proof of (3.3) is finished. The proof of (3.4) and (3.5) are completely analogous and 
will be left to the reader. $\blacksquare$
\begin{coro} Let $(M^{2n},g,J)$ be a compact K\"ahler manifold with a $2$-formal Riemannian metric. Then : \\
(i) The length of any harmonic $2$-form is constant over $M$; \\
(ii)  If $\alpha$ in ${\cal{H}}_{0}^{1,1}$ has vanishing square, then $\alpha$ is necessarily $0$.
\end{coro}
{\bf{Proof}} : \\
Indeed, if $\alpha$
belongs to ${\cal{H}}^{1,1}$ or to ${\cal{H}}_{-}^2$ then we saw that $S^2$ belongs to ${\cal{S}}$ and all elements of ${\cal{S}}$ have constant 
trace. But the trace of $S^2$ equals
the squared norm of $\alpha$. The other statement is straightforward.\\
$\blacksquare$
\begin{pro}
Suppose that $(M^{2n},g,J)$ is a compact K\"ahler manifold such that the metric $g$ is $2$-formal. Let $\alpha$ belong to ${\cal{H}}^{1,1}$ and let
$S$ in ${\cal{S}}$ be the associated symmetric endomorphism of $TM$. Then : \\
(i) The eigenvalues of $S$ are constant with eigenbundles of constant rank. \\
(ii) If $\lambda_i, 1 \le i \le p$ are (the pairwise distinct ) eigenvalues of $S$ we have an orthogonal and $J$-invariant 
decomposition
$$ TM=\bigoplus \limits_{j=1}^{p}E_i $$
where $E_i$ is the eigenspace of $S$ corresponding to $\lambda_i$. Furthermore, for all $1 \le i \le p$ the distributions
$E_i$ and $\hat{E}_i=\bigoplus \limits_{j=1, j \neq i}^{p}E_i$ are integrable. 
\end{pro}
{\bf{Proof}} : \\
(i) From (3.3) we deduce that $S^k$ belongs to ${\cal{S}}$ for all $k$ in $\mathbb{N}$. As ${\cal{S}}$ is finite
dimensional there exists $P$ in $\mathbb{R}[X]$ such that $P(S)=0$. Since $S$ is symmetric, $P$ can be supposed to 
have only real roots and again by the symmetry of $S$ we can moreover assume that all these roots are simple. Let $\lambda_i, 1 \le i \le p$ be
these (pairwise distinct) roots and let $m_i$ be the dimension of the corresponding eigenbundle. To see that 
$m_i, 1 \le i \le p$ are constant over $M$ we use the fact that $S^k$ belongs to ${\cal{S}}$ for all $k$ in $\mathbb{N}$ in order to deduce that $Tr(S^k)=c_k$ for some constant $c_k$ and
for all natural $k$. In other words
$$ \sum \limits_{i=1}^{p}m_i \lambda_i^k=c_k$$
for all $k$ in $\mathbb{N}$. Solving this Vandermonde system leads to the constancy of the functions $m_i, 1 \le
i \le p$.\\
(ii) Let $\omega^i$ be the orthogonal projection of the K\"ahler form $\omega$ on $E_i, 1 \le i \le p$. Then $\alpha=
\sum \limits_{i=1}^{p} \lambda_i \omega^i$ and moreover, by (3.3) we obtain that 
$$  \sum \limits_{i=1}^{p}\lambda_i^k\omega^i$$
belongs to ${\cal{H}}^{1,1}$ for all natural $k$. We assume now that the eigenvalues of $S$ are ordered by  
$ \vert \lambda_1 \vert  < \vert  \lambda_2 \vert <  \ldots < \vert \lambda_p \vert $. We divide by
$ \vert \lambda_p \vert^k $ and make $k \to \infty$. It follows that $\omega_p$ belongs to ${\cal{H}}^{1,1}$. By 
induction the same holds for $\omega^i, 2 \le i \le p$. If $\lambda_1 \neq 0$, the form $\omega^1$ is trivially in
${\cal{H}}^{1,1}$ and if $\lambda_1=0$ the same is true since $\omega=\sum \limits_{i=1}^{2n} \omega^i$. \par
Therefore $\omega^i, 1 \le i \le p$ are all closed. Fix $1 \le i \le p$ and consider the decomposition $TM=E_i \oplus
\hat{E}_i$. Let $X,Y$ be in $E_i$ and $V$ in $\hat{E}_i$. A straightforward computation  yields to 
$(\nabla_X\omega^i)(Y,V)=-<\nabla_XY, JV>, (\nabla_Y\omega^i)(X,V)=-<\nabla_YX, JV>$ and
$(\nabla_V \omega^i)(X,Y)=0$. Now the closedeness of $\omega^i$ ensures the integrability of $E_i$. That of 
$\hat{E_i}$ is proved in a similar way, by computing $(d\omega^i)(V,W,X)$ with $W$ in $\hat{E}_i$.
$\blacksquare$ 
$\\$
\par
As an immediate consequence of Proposition 3.1 we obtain : 
\begin{coro}
Let $(M^{2n},g,J)$ be a compact K\"ahler manifold such that $g$ is $2$-formal. Then any element of ${\cal{H}}^{1,1}$ can be uniquely written as a linear
combination of $g$-compatible symplectic forms.
\end{coro}
{\bf{Proof}} : \\
Proposition 3.1 actually says that any $2$-form $\alpha$ in ${\cal{H}}^{1,1}$ can be written as 
$\alpha=\sum \limits_{j=1}^{p}\lambda_i \omega^i$ where $\omega^i$ denotes the projection of the K\"ahler form
$\omega$ on $E_i, 1 \le i \le p$. Moreover the forms $\omega^i$ belong to ${\cal{H}}^{1,1}$ for all $1 \le i \le p$. Now, for $1 \le k \le p$ 
we define an almost complex
structure $J_k$ on $TM$ by setting 
$$ J_k=J \ \mbox{on} \ \hat{E}_k,  \ \mbox{and} \ J_k=-J \ \mbox{on} \ E_k. $$
An easy consequence of the integrability of {\it{both}} $E_k$ and $\hat{E}_k$ is that $(g, J_k)$ are almost K\"ahler structures 
(i.e. the corresponding K\"ahler forms 
$\Omega_k=g(J_k \cdot, \cdot)$ are closed) commuting with $J$, for any
$1 \le k \le p$. Note that $J_k$ is integrable, i.e. $(g, J_k)$ is a K\"ahler structure iff $E_k$ is parallel with respect 
to the Levi-Civita connection. To finish the proof of the Corollary it suffices to note that 
$$ \omega^k=\frac{1}{2} \omega-\frac{1}{2}\Omega_k$$
for all $1 \le k \le p$.
$\blacksquare$ \par
$\\$
Therefore the proof of part (i) of Theorem 1.2 is now complete.

\section{More on Hodge numbers}
The aim of this section is to provide some information about the Hodge numbers of a geometrically formal K\"ahler manifold. We begin with the following simple
observation.
\begin{pro}
Let $(M^{2n},g, J)$ be a compact K\"ahler manifold such that the metric $g$ is formal. Then $h^{0,n}(M) \le 1$ and equality holds iff $g$ is a Ricci flat metric.
\end{pro}
{\bf{Proof}} : \\
Because any harmonic is of constant length the Hodge numbers of $(M^{2n},g,J)$ are bounded by the dimensions of their corresponding vector bundles 
hence $h^{0,n}(M) \le 1$. If equality holds, it follows
that the canonical bundle of $(M^{2n},g,J)$ is trivialized by a harmonic $(0,n)$-form of constant length and this leads in the standard way to 
the vanishing of the Ricci tensor.
$\blacksquare$ \\
$\\$
We investigate now the structure of $J$-anti-invariant harmonic $2$-forms.
\begin{pro}
Let $(M^{2n},g,J)$ be a compact K\"ahler manifold such that the metric $g$ is $2$-formal. Then  : \\
(i) Any non-zero element of ${\cal{H}}_{-}^{2}$ induces
in a canonical way a local splitting of $M$ as the Riemannian product of a compact K\"ahler manifold $M_1$ and
a compact hyperk\"ahler manifold $M_2$. \\
(ii) ${\cal{H}}_{-}^2$ consists only in parallel forms. 
\end{pro}
{\bf{Proof}} : \\
(i) Let $\alpha$ be in ${\cal{H}}^2_{-}$ be non-zero and let $A$ in ${\cal{A}}$ be its associated endomorphism. Then $A^2$
belongs to ${\cal{S}}$ by (3.4) and using proposition 3.1 we obtain a $J$-invariant and orthogonal decomposition 
$$ TM=\bigoplus \limits_{i=1}^{p} E_i$$
where $E_i$ are eigenspaces of $A^2$ for the (constant) eigenvalues 
$\mu_i \le 0, 1 \le i \le p$. Now using again (3.4) we get
that $A^{2k}$ belongs to ${\cal{S}}$ and further, by (3.5) that $A^{2k+1}$ is in ${\cal{A}}$ for all $k$ in $\mathbb{N}$. Let 
$A_i, 1 \le i \le p$ be the orthogonal projections of $A$ on $E_i$. An argument similar to the proof of
Proposition 3.1, (i) shows that $A_i$ are in ${\cal{A}}$ for all $1 \le i \le p$. Equivalently, the forms $\alpha^i$, 
associated to $A_i$ are in ${\cal{H}}_{-}^2$ for all $1 \le i \le p$ and therefore have to be closed. We will show now
that one can reduces to the case when  $A$ has no kernel. Indeed, let us assume that $A$ has non-empty
kernel, that is $A^2$ has a zero eigenvalue, say $\mu_1$. Set ${\cal{V}}=E_1$ and $H=\hat{E}_1$. Then the
endomorphism
$F=\sum \limits_{i=2}^{p}\frac{1}{\sqrt{-\mu_i}}A^i$, an element of ${\cal{A}}$, vanishes on ${\cal{V}}$ and 
defines an almost complex structure $I$ on $H$, compatible with $g$ and such that $IJ+JI=0$. Since the $2$-form
form associated to $F$ is closed we get :
\begin{nr} \hfill
$ <(\nabla_{U_1}F)U_2, U_3>-<(\nabla_{U_2}F)U_1, U_3>+<(\nabla_{U_3}F)U_1, U_2>=0\hfill $
\end{nr}
for all $U_j, 1 \le j \le 3$ in $TM$. Let $\overline{\nabla}$ be the metric connection leaving ${\cal{V}}$ and $H$ 
parallel. Then $\nabla_XY=\overline{\nabla}_XY+A_XY$ for all $X,Y$ in $H$ where the O'Neill-type tensor
$A : H \times H \to {\cal{V}}$ is the obstruction to the distribution $H$ to be totally geodesic. Taking
$U_1=X, U_2=Y$ and $U_3=V$ in (4.1) with $X,Y$ in $H$ and $V$ in ${\cal{V}}$ we get :
\begin{nr} \hfill
$ <A_{X}IY-A_Y(IX),V>+<(\overline{\nabla}_VI)X,Y>=0.\hfill $ 
\end{nr}
Since the connection $\overline{\nabla}$ is metric and $I^2=-1$ on $H$ it follows that $<(\overline{\nabla}_VI)IX,IY>=
-<(\overline{\nabla}_VI)X,Y>$. Therefore, changing $X$ in $IX$ and $Y$ in $IY$ in (4.2) and summing the result with
(4.2) we obtain $A_{X}IY-A_Y(IX)-A_{IX}(Y)+A_{IY}X=0$. But $A$ is symmetric as $H$ is 
integrable (see Proposition 3.1, (ii)) hence
\begin{nr} \hfill
$A_X(IY)=A_Y(IX) \hfill $
\end{nr} for all $X,Y$ in $H$. As $(g,J)$ is K\"ahler and, by construction both
${\cal{V}}$ and $H$ are $J$-invariant, we are lead to $A_X(JY)=JA_XY$ for all $X,Y$ in $H$. Taking this into
account and replacing $Y$ by $JY$ in (4.3) yields after
a standard manipulation to the vanishing of $A$. \par
We showed that $H$ is a totally geodesic distribution hence the foliation induced by ${\cal{V}}$ is
a Riemannian one. Now, 
on any integral manifold of $H$, with respect of the induced metric, the triple
$I,J,K=IJ$ induce a family of almost complex structures satisfying the quaternionic identities and with
closed associated K\"ahler forms. Then a well
known lemma due to Hitchin \cite{Hit} implies that the metric is hyperk\"ahler and hence Ricci flat. It follows that the
transversal Ricci curvature $Ric^H$ of the Riemannian foliation induced by ${\cal{V}}$ vanishes and using
Proposition 2.1 we obtain that ${\cal{V}}$ is also totally geodesic, hence the desired splitting. \\
(ii) It suffices to work on the compact, Ricci flat manifold $M_2$ where $A$ has no kernel. Then any of the commuting
almost K\"ahler structures induced by $A^2$ have to be K\"ahler by a theorem of Sekigawa \cite{seki1}. This
means that the spaces $E_i $ are all parallel and on each of them $A$ is a multiple of a K\"ahler
structure (anti-commuting with $J$). This implies immediately the parallelism of $\alpha$.
$\blacksquare$  \\ 
$\\$
In particular Theorem 1.2 is now completely proved. It can be used to refine, in the K\"ahler case, the Betti number estimates 
$b_p(N) \le b_p(T^n)$ known to hold (see \cite{Kot1}) for an arbitrary geometrically formal manifold $(N^n,h)$.
\begin{coro}Let $(M^{2n},g,J)$ be K\"ahler such the metric $g$ is formal. If $M$ is locally irreducible
 then $b_1(M)=0$ and $b_{2p+1}(M) \le C_{2n}^{2p+1}-2n$ for all $p \ge 1$.
\end{coro}
{\bf{Proof}} : \\
This is an immediate consequence of theorem 1.1 and of the Lefschetz decomposition (see \cite{Gold}) of the harmonic
forms of a K\"ahler manifold .
$\blacksquare$ \\ $\\$ 
The previous corollary can also be reformulated to say that if a geometrically formal K\"ahler manifold has a maximal
Betti number of odd degree then the metric is a flat one. More generally we have : 
\begin{coro}Let $(M^{2n},g,J)$ be K\"ahler such that the metric $g$ is formal. If there exists $1 \le p \le 2n-1$ such that
$b_p(M)=b_p(\mathbb{T}^{2n})$ then $g$ is flat metric.
\end{coro}
{\bf{Proof}} : \\
By Corollary 4.1 it suffices to study the case $p=2q$. In view of the Hodge duality we may also suppose that
$p \le n$. Using the fact that harmonic forms of $g$ are of constant length and the Lefschetz decomposition of a
K\"ahler manifold, we are lead to $b_2(M)=b_2(\mathbb{T}^n)$ and further to $b_2^{-}(M)=
b_2^{-}(\mathbb{T}^n)$. But in the case of the torus it is an algebraic fact that $\{{\cal{A}}, {\cal{A}} \}={\cal{S}}$ hence
in view of the parallelism of $J$-anti-invariant harmonic $2$-forms ${\cal{H}}^{1,1}$ equally consists of parallel
forms. We have therefore a framing of $\Lambda^2(M)$ by parallel
two-forms and this implies in a standard way the desired result.
$\blacksquare$ \\ \par
The rest of the section will be consecrated to explore a number of consequences of Theorem 1.2 under
various curvature assumptions. First of all we have :
\begin{coro}
Suppose that $(M^{2n},g,J)$ is a compact quotient of the complex hyperbolic space, endowed with its canonical
K\"ahler metric. If $g$ is a formal metric then $b_1(M)=0$ and $b_2(M)=1$.
\end{coro}
{\bf{Proof}} : \\
Since $(M^{2n},g,J)$ is locally irreducible, it follows that $b_1(M)=0$ by Theorem 1.1 and also that $b_2^{-}(M)=0$ by Proposition 4.2. Now using Corollary 3.2 and
the fact (see \cite{Nagy}) that on $(M^{2n},g,J)$ every orthogonal, $J$-commuting, almost K\"ahler structure has to be K\"ahler and therefore a multiple
of $J$ we are lead to $h^{1,1}=1$.
$\blacksquare$ \\ \par
We investigate now the incidence of having constant scalar curvature on $2$-formal K\"ahler metrics.
\begin{teo}
Let $(M^{2n},g,J)$ be a compact K\"ahler manifold. Assume that the metric $g$ is $2$-formal and locally
irreducible. Then : \\
(i) if the scalar curvature of $g$ is constant then the eigenvalues of $Ric$ are constant (together with their
multiplicities) over $M$. \\
(ii) If the scalar curvature is constant and $Ric_g \ge 0$ then $h^{1,1}=1$. Moreover, under these assumptions $g$ is an Einstein metric.
\end{teo}
{\bf{Proof}} : \\
(i) If the scalar curvature is constant, the Ricci form is harmonic and the result follows by Proposition 3.1, (i). \\
(ii) In this case the Ricci tensor has only constant and non-negative eigenvalues, by (i). Using Proposition 3.1, (ii) we can
always rescale, by an argument similar to Lemma 2.2, page 774, in \cite{Apo1}, the metric along the eigenbundles of $Ric$
in order to get a K\"ahler metric with $2$ constant and non-negative eigenvalues. Then the splitting result of \cite{Apo1}
asserts that every $g$-compatible almost K\"ahler structure, commuting with $J$ is in fact K\"ahler. Therefore,
by local irreducibility, $h^{1,1}=1$ leading further, by (i) to the fact that $g$ is Einstein.
$\blacksquare$ \\ $\\$
An immediate consequence of Proposition 4.2 and Theorem 4.1, (ii) is the following.
\begin{teo} Let $(M^{2n},g,J)$ be a compact K\"ahler manifold such that $g$ is $2$-formal and locally
irreducible. Then either  \\
(i) $b_2^{-}(M)=0$ \\
or \\
(ii) $b_2^{-}(M)=2$ and $(M,g)$ is a hyperk\"ahler manifold. Moreover, in this case we must have $h^{1,1}=1$ and
thus $b_2(M)=3$.
\end{teo}
Thus, a naturally arising question is to decide whether a hyperk\"ahler metric can be geometrically formal. This
seems quite unlikely from the perspective that known examples of compact hyperk\"ahler manifolds
(see \cite{Joyce} for an account) have second Betti number greater than $3$. However, we were unable to prove that
hyperk\"ahler metrics cannot be geometrically formal, except in lows dimensions, as the following shows :
\begin{pro}
They are no geometrically formal and locally irreducible hyperk\"ahler manifolds in dimensions $4$ and $8$.
\end{pro}
{\bf{Proof}} : \\
In dimension $4$ this was proven in \cite{Kot}. To prove the statement in dimension $8$ we need to recall some
facts about the topology of hyperk\"ahler manifolds. Thus, let $Z^{4m}$ be a hyperK\"ahler manifold. It was proven in
\cite{Sal1} that the Betti numbers of $Z$ satisfy the following remarkable relation
\begin{nr} \hfill
$ 3P^{\prime \prime}(-1)=m(12m-5)P(-1)  \hfill $
\end{nr}
where $P(t)=\sum \limits_{k=0}^{4m}b_{k}(Z)t^{k}$ is the Poincar\'e polynomial. \par
Suppose now that $(M^8, g, I, J, K)$ is a hyperk\"ahler manifold such that the metric $g$ is formal. Then
$b_1(M)=0$ and $b_2(M)=3$ hence after an easy computation (4.4) becomes
$$ b_3(M)+b_4(M)=76.$$
To obtain a contradiction we will produce estimates of the Betti numbers $b_3$ and $b_4$.
We denote by $\omega_I, \omega_J, \omega_K$ the
corresponding K\"ahler forms. Since $b_1(M)=0$, any harmonic $3$-form lives in the orthogonal
complement of $\{ \alpha \wedge \omega_I+\beta \wedge \omega_J +
\gamma \wedge \omega_K : \alpha, \beta, \gamma \in \Lambda ^1(M)\}$. Since any harmonic $3$-form has constant length
we get $b_3(M) \le C_8^3-3\cdot 8=32$. Let us denote by $\Lambda_{o}^2(M)$ the subbundle of $\Lambda ^2(M)$
consisting of forms orthogonal to
$\omega_I, \omega_J, \omega_K$. This bundle do not contain any harmonic form (since the second
cohomology group is generated by the hyperk\"ahler forms) and therefore any harmonic $4$-form must be orthogonal to
$E$, the subbundle of $\Lambda^4(M)$ generated by $\{ \omega_I \wedge \alpha+\beta \wedge \omega_J +\gamma \wedge \omega_K : \alpha,
\beta, \gamma \in \Lambda ^2_{o}(M)\}$. Now, a simple computation shows that $\omega_I \wedge \alpha$ and $\omega_J \wedge \beta$ have to be
orthogonal for all $\beta$ in $\Lambda^2_{o}(M)$ and
$\alpha$ in $\Lambda^{1,1}_{K}(M) \cap \Lambda^2_{o}(M)$. Therefore the rank of $E$ is greater than
$(C_8^2-3)+15=40$ and this implies, again by the fact that
harmonic $4$-forms are of constant length, that $b_4(M) \le C_8^4-40=30$. We found that
$b_3(M) +b_4(M) \le 62$, an obvious impossibility.
$\blacksquare$ \\ \par

\section{$6$-dimensions}
We consider in this section a $6$-dimensional K\"ahler manifold $(M^6,g,J)$
such that the metric $g$ is formal and locally irreducible. Our aim is to give the possible values for
the second Betti number and also to explore the algebraic structure of the
second cohomology group. This will also lead to the proof of Theorem 1.3. We first prove : 
\begin{lema} Let $(M^6,g,J)$ be a geometrically formal, K\"ahler
manifold. If $g$ is locally irreducible then : \\
(i) $b_1(M)=b_2^{-}(M)=0$. \\
(ii) we have either $Td(M)=1$ or $Td(M)=0$. If the last case occurs then $g$ is Ricci flat and $b_2(M)=1$.
\end{lema} 
{\bf{Proof}} : \\ 
(i) Direct consequence of Theorem 1.1 and Proposition 4.2, (i). \\
(ii)  As $h^{0,2}(M)=0$ the use of Riemann-Roch tells us that $Td(M)=1-h^{0,3}$. Using Proposition 4.1 we get that
either $h^{0,3}(M)=0$ (and hence $Td(M)=1$) or $h^{0,3}(M)=1$ (thus $Td(M)=0$) and $g$ is Ricci flat. But if $g$ is Ricci flat
one uses Theorem 4.1, (ii) to get that $b_2(M)=1$ and th proof is finished.
$\blacksquare$ \\ $\\$
We shall now study commutation rules inside  ${\cal{H}}^{1,1}$. Our methods will be mainly topological and shall rely on the following : 
\begin{pro} Let $(M^6,g,J)$ be a compact K\"ahler
manifold. Suppose that we have an orthogonal and $J$-invariant decomposition
$TM={\cal{V}} \oplus H$ where ${\cal{V}}$ is of real rank $2$. If $I_1$ and
$I_1$ are almost complex structures on $H$ which are compatible with $g$ and
such that $\{ I_1, I_2\}=0$ and $[I_k, J]=0, k=1,2$ then $\chi(M)$ is
divisible by $12$.
\end{pro}
{\bf{Proof}} : \\ We have
$24Td(M,J)=c_1(M,J)c_2(M,J)$. Or $c_1(M,J)=c_1({\cal{V}})+c_1(H)$ and
$c_2(M,J)=c_2(H)+c_1({\cal{V}})c_1(H)$ where the bundles ${\cal{V}}$ and $H$
are endowed with the complex structure induced by $J$. Then :
\begin{nr} \hfill $
24Td(M,J)=c_1({\cal{V}})c_2(H)+c_1(H)c_2(H)+c_1^2({\cal{V}})c_1(H)+c_1({\cal{V}})c_1^2(H). \hfill $ 
\end{nr} 
Consider now the almost complex structure $J_0$ which equals $-J$ on ${\cal{V}}$ and $J$ on $H$. Taking into account that $J_0$ is 
inducing the orientation opposite to that induced by $J$ we get as before : 
$$ -24Td(M,J_0)=-c_1({\cal{V}})c_2(H)+c_1(H)c_2(H)+c_1^2({\cal{V}})c_1(H)-c_1({\cal{V}})c_1^2(H).$$ Subtracting we obtain : 
$$ 12(Td(M,J)+Td(M,J_0))=c_1({\cal{V}})(c_2(H) +c_1^2(H)).$$ 
Consider now the orthogonal involution $\sigma=JI_1$ of $H$. It can be used to obtain an orthogonal and $J$-invariant decomposition 
$H=H^{+} \oplus H^{-}$ where
$H^{\pm}$ are the $\pm 1$-eigenspaces of $\sigma$. Since $\{\sigma, I_2 \}=0$
we have that $H^{-}=I_2H^{+}$ hence $I_2$ induces a complex isomorphism
between $H^{+}$ and $H^{-}$. It follows that $c_1(H^{+})=c_1(H^{-})$ and
therefore $c_1(H)=2c_1(H^{+})$ and $c_2(H)=c_1^2(H^{+})$. We deduce that
$c_1^2(H)=4c_2(H)$ and further : $$
12(Td(M,J)+Td(M,J_0))=5c_1({\cal{V}})c_2(H).$$ Now
$\chi(M)=c_3(M,J)=c_1({\cal{V}})c_2(H)$ and this leads to
$12(Td(M,J)+Td(M,J_0))=5\chi(M)$, finishing the proof of the proposition. 
$\blacksquare$ \\ \par
In order to prove Theorem 1.3 we need a number of preliminary results. 
\begin{lema}Let
$(M^6,g,J)$ be a geometrically formal K\"ahler manifold and let $\sigma$ in
${\cal{S}}$ be an involution. Then : \\ (i) We have ${\cal{S}}=C_{\sigma}
\oplus A_{\sigma}$ where we defined  $ C_{\sigma}=\{ S \in {\cal{S}} : [ S,
\sigma ]=0 \}$ and $ A_{\sigma}=\{ S \in {\cal{S}} : \{ S, \sigma \}=0 \}$. \\
(ii) the dimension of $A_{\sigma}$ is less or equal to $1$. 
\end{lema} 
{\bf{Proof}} : \\
(i) Let $S$ be in ${\cal{S}}$ and decompose $S=S_1+S_2$ where $S_1$ and $S_2$ are commuting resp. 
anti-commuting with $\sigma$. To prove the result it is enough to see that $S_1$ belongs to ${\cal{S}}$. But using
(3.3) we obtain that $S\sigma+\sigma S$ is in ${\cal{S}}$ and again by (3.3), 
$\{\sigma, S\sigma+\sigma S \}=2(S+\sigma S \sigma  )=4S_1$ belongs to ${\cal{S}}$ and the proof is finished. \\
(ii) Assume that $g$ is not Ricci flat because otherwise $b_2(M)=1$ (see Lemma 5.1, (ii)) and there is nothing to prove. 
Let $TM={\cal{V}} \oplus H$ be the 
orthogonal and $J$-invariant decomposition of $TM$ in the $-1$ and $1$-eigenspaces of $\sigma$, and let us
assume that ${\cal{V}}$ has real rank $2$. Suppose that $A_{\sigma}$ is non-empty and let 
$S$ be a non-vanishing element of $A_{\sigma}$. Then
$S({\cal{V}}) \subseteq H$ and hence $S$ defines an $J$-invariant isomorphism from ${\cal{V}}$ to its  
image $H_1$. If $S^{\prime}$ in $A_{\sigma}$ is orthogonal to $S$ then it defines a $J$-invariant isomorphism
from ${\cal{V}}$ to $H_2$ ,the orthogonal complement of $H_1$ in $H$. In other words we have a decomposition
$TM={\cal{V}} \oplus H_1 \oplus H_2$ in $J$-isomorphic bundles. Let us denote by $h$ the first Chern class of
${\cal{V}}$. Of course, $c_1(H_1)=c_1(H_2)=h$. From (5.1) we deduce easily that $24Td(M,J)=9h^3$ and moreover
$\chi(M)=c_3(M)=h^3$. But using Lemma 5.1, (ii) we infer that $Td(M)=1$ and this leads to $24=9\chi(M)$ and since this equation
has no integer solution we obtained a contradiction with the existence of $S^{\prime}$, hence finishing the proof of the lemma.
$\blacksquare$
\begin{rema}
From the above proof, we see that Lemma 5.2 continues to hold for $2$-formal metrics provided
that the Todd genus $Td(M,J)$ is not divisible by $3$.
\end{rema}
Now we need another auxiliary result in order to get some precisions concerning the structure of 
harmonic $3$-forms in some special cases.  By contrast with the previous remark, from now on we will use the formality
hypothesis in a crucial way.
\begin{lema}
(i) Let $(M^6,g)$ be geometrically formal and let $\alpha$ be a harmonic $2$-form. If $L_{\alpha}$ is the exterior multiplication with $\alpha$ then its 
adjoint, to be denoted by $L_{\alpha}^{\star}$ preserves the space of harmonic forms. \\
(ii) Let $(M^6,g,J)$ be an almost K\"ahler manifold, which is also geometrically formal. The for all primitive $\alpha$ in ${\cal{H}}^3$ we have that 
$J\alpha$ is also harmonic. Here $J$ acts on a $3$-form $\beta$ by $(J\beta)(X,Y,Z)=\beta(JX,JY,JZ)$ whenever $X,Y,Z$ are in $TM$.
\end{lema}
{\bf{Proof}} : \\
(i) Using the definition of the Hodge star operator we see that $L_{\alpha}^{\star}$ is up to a sign equal to the composition $\star L_{\alpha} \star$ and the
result follows by the formality hypothesis. \\
(ii) It suffices to note, using for instance the definition of the Hodge-star operator, that $\star \alpha=J\alpha$ and the result follows.
$\blacksquare$ \\ $\\$
Our last preparatory result consists in giving a bound on the third Betti number in case that, for some 
involution $\sigma$ of ${\cal{S}}$ the space $A_{\sigma}$ is non-empty.
\begin{lema}
Let $(M^6,g,J)$ be a geometrically formal K\"ahler manifold and let $\sigma$ in ${\cal{S}}$ be an involution. If $ A_{\sigma}$ is $1$-dimensional then : \\
(i) $b_3(M) \le 6$ \\
(ii) $12 \vert \chi(M)$.
\end{lema}
{\bf{Proof}} : \\
(i) Let $TM={\cal{V}} \oplus H$ be the spectral decomposition of $\sigma$ into $\pm$-eigenbundles and suppose furthermore that
the real rank of ${\cal{V}}$ equals $2$. It is easy to verify that we are then in the situation of Proposition 5.1, i.e. we have almost complex structures
$I_1, I_2$ on $H$ which are mutually anti-commuting and commuting with $J$. Let us consider now $\alpha$ a harmonic $3$-form. If $\omega_1$ is the orthogonal
projection of $\omega$ on ${\cal{V}}$ we have that $L^{\star}_{\omega_1} \alpha=0$ hence we can decompose $\alpha=\alpha_1+\alpha_2$ where
$\alpha_1 $ belongs to $\Lambda^1({\cal{V}}) \otimes \Lambda^2(H)$ and $\alpha_2$ is in $\Lambda^3(H)$. Then $J\alpha=J\alpha_1+J\alpha_2$ belongs to
${\cal{H}}^3$ and also $J_1 \alpha=J_1 \alpha_1+J_1 \alpha_2=-J\alpha_1+J\alpha_2$ belongs to ${\cal{H}}^3$ where $J_1$ is the almost K\"ahler structure on $M$ acting as
$-J$ on ${\cal{V}}$ and as $J$ on $H$. It follows that $J\alpha_2$ is in $H$ and thus $\alpha_2=-J(J\alpha_2)$ is a harmonic $3$-form. If $\omega_2$ is the orthogonal
projection of $\omega$ on $H$ then $L^{\star}_{\omega_2}\alpha_2=0$ meaning that $\alpha_2$ vanishes. \par
We showed that any harmonic $3$-form $\alpha$ belongs to $\Lambda^1({\cal{V}}) \otimes \Lambda^2(H)$ and in order to prove the lemma, we have to take into account the
hypothesis that $A_{\sigma}$ is $1$-dimensional. Look at $\alpha$ as a vector bundle morphism $\alpha : \Lambda^2(H) \to {\cal{V}}$. Since $L^{\star}_{\omega_2}\alpha=
L^{\star}_{\omega_{I_1}}\alpha=L^{\star}_{\omega_{I_2}}\alpha=0$ we find that $\alpha$ is in fact defined from $\Lambda_{-}^{2}(H) \oplus E$ to ${\cal{V}}$ where
$E$ is the orthogonal complement of the span of $\omega_{I_1}, \omega_{I_2}$ in $\Lambda_0^{1,1}(H)$. But all harmonic $3$-forms are of constant length hence the third Betti
number cannot exceed $6$, the rank of the vector bundle $(\Lambda_{-}^{2}(H) \oplus E) \otimes {\cal{V}}$.\\
(ii) follows immediately from Proposition 5.1.
$\blacksquare$ \\ $\\$
We are now in position to give the \\
$\\$
{\bf{Proof of Theorem 1.3}} : \\
By Lemma 5.1 and Corollary 3.2 we only need to prove the bound on $b_2(M)$. We can suppose that $b_2(M) \ge 2$, otherwise there is nothing to prove. It follows 
that ${\cal{S}}$ contains a non-trivial involution $\sigma$ of $TM$ whose 
$\pm 1$-eigenspaces will be denoted by ${\cal{V}}$ and $H$. Moreover we can suppose that ${\cal{V}}$ has real
rank $2$. Let us suppose now that $dim_{\mathbb{R}}A_{\sigma}=1$
and let $S$ in $A_{\sigma}$ be non-zero. Then we have $TM={\cal{V}} \oplus H_1 \oplus H_2$ with 
$H_1=S({\cal{V}})$ and $H_2$ the orthogonal complement of $H_1$ in $H$. Let now $S_1$ be in $C_{\sigma}$. Then $SS_1+S_1S$ is in 
$A_{\sigma}$ and then $SS_1+S_1S=\lambda S$ where $\lambda$ is a real constant. It follows by some simple algebraic considerations that
$S_1$ preserves $H_i, i=1,2$ and since the latter are $2$-dimensional we obtain that $C_{\sigma}$ has dimension $3$ hence $b_2(M)=4$ by Lemma 5.2, (i). By 
Lemma 5.4, (i) and the
fact that $(M^6,g,J)$ is K\"ahler, the only possibilities for $b_3(M)$ are $0,2,4,6$. Therefore, the possible values of 
$\chi(M)=2+2b_2(M)-b_3(M)=10-b_3(M)$ are $10, 8, 6, 4$ a fact which in contradiction with the fact that $\chi(M)$ is divisible by $12$ (cf. Lemma 5.4, (ii)). \par
We showed that for any involution in ${\cal{S}}$ the space $A_{\sigma}$ vanishes and this implies that any two elements of ${\cal{S}}$ must commute. At a given
point of $M$, the elements of ${\cal{H}}^{1,1}$ form a commutative subalgebra of $\mathfrak{u}(3)$ and this implies finally $b_2(M) \le 3$.
$\blacksquare$ \\
\begin{rema}
(i) In view of the Theorem 1.3 it seems
quite likely that a case by case discussion could give the real cohomology type of a geometrically formal, $6$-dimensional K\"ahler manifold provided that one founds
a method to analyze obstructions to geometric formality at the level of the third cohomology group.  At this moment all the informations about the third Betti
number we have are the estimates $b_3(M) \le 10$ if $b_2(M)=2$ and $b_3(M) \le 8$ if $b_3(M)=3$; these follow easily from the first part of the proof of lemma 5.4.\\
(ii) If $M$ is K\"ahler and geometrically formal of dimension divisible  by $4$, the commutativity result of Theorem 1.3 may not hold since, a priori, $(g,J)$  could
admit a compatible, complex symplectic structure, which is also $J$-invariant.
\end{rema} 

{\bf{Acknowledgement}} : The author wishes to thank the referee for his helpful comments.

\center
\begin{flushright}
Paul-Andi Nagy \\
Humboldt Universit\"at zu Berlin \\
Institut f\"ur Mathematik\\
Sitz : Rudower Chaussee 25, D-10099 Berlin.\\
email : nagy@mathematik.hu-berlin.de
\end{flushright}
\end{document}